\documentclass[11pt,leqno]{amsart}
\usepackage[english]{babel}
\usepackage{hyperref}
\usepackage{amssymb}
\usepackage[mathcal]{eucal}

\numberwithin{equation}{section}

\theoremstyle{definition}
\newtheorem{ex}{Example}[section]

\theoremstyle{plain}
\newtheorem{prop}[ex]{Proposition}
\newtheorem{teor}[ex]{Theorem}
\newtheorem{cor}[ex]{Corollary}

\theoremstyle{remark}
\newtheorem*{oss-tech}{Remark}
\newtheorem*{org-not}{Organization and Notation}

\newcommand{\N}{\mathbb{N}}
\newcommand{\Z}{\mathbb{Z}}

\newcommand{\F}{\mathbb{F}}

\DeclareMathOperator{\rk}{rk}
\DeclareMathOperator{\dd}{d}
\DeclareMathOperator{\Zent}{Z}
\DeclareMathOperator{\Cen}{C}

\DeclareMathOperator{\Aut}{Aut}
\DeclareMathOperator{\Kern}{ker}

\begin{document}

\title[Finite axiomatizability of rank and dimension]{Finite
  axiomatizability of the rank and \\ the dimension of a pro-$\pi$
  group}

\author[M. Conte]{Martina Conte} 
\address{Martina Conte:
  Mathematisches Institut, Heinrich-Heine-Universit\"at, 40225
  D\"usseldorf, Germany} \email{conte@uni-duesseldorf.de}

\author[B. Klopsch]{Benjamin Klopsch} 
\address{Benjamin Klopsch:
  Mathematisches Institut, Heinrich-Heine-Universit\"at, 40225
  D\"usseldorf, Germany} \email{klopsch@math.uni-duesseldorf.de}



\dedicatory{In Memory of Avinoam Mann}

\thanks{The research was conducted in the framework of the
  DFG-funded research training group ``GRK 2240: Algebro-Geometric
  Methods in Algebra, Arithmetic and Topology''.}


\begin{abstract}
  The Pr\"ufer rank $\rk(G)$ of a profinite group $G$ is the
  supremum, across all open subgroups $H$ of $G$, of the minimal
  number of generators $\dd(H)$.  It is known that, for any given prime
  $p$, a profinite group $G$ admits the structure of a $p$-adic
  analytic group if and only if $G$ is virtually a pro-$p$ group of
  finite rank.  The dimension $\dim G$ of a $p$-adic analytic
  profinite group $G$ is the analytic dimension of $G$ as a $p$-adic
  manifold; it is known that $\dim G$ coincides with the rank $\rk(U)$
  of any uniformly powerful open pro-$p$ subgroup $U$ of $G$.

  Let $\pi$ be a finite set of primes, let $r \in \N$ and let
  $\mathbf{r} = (r_p)_{p \in \pi}, \mathbf{d} = (d_p)_{p \in \pi}$ be
  tuples in $\{0, 1, \ldots,r\}$.  We show that there is a single
  sentence $\sigma_{\pi,r,\mathbf{r},\mathbf{d}}$ in the first-order
  language of groups such that for every pro-$\pi$ group $G$ the
  following are equivalent: (i) $\sigma_{\pi,r,\mathbf{r},\mathbf{d}}$
  holds true in the group $G$, that is,
  $G \models \sigma_{\pi,r,\mathbf{r},\mathbf{d}}$; (ii) $G$ has rank
  $r$ and, for each $p \in \pi$, the Sylow pro-$p$ subgroups of $G$
  have rank $r_p$ and dimension $d_p$.

  Loosely speaking, this shows that, for a pro-$\pi$ group $G$ of
  bounded rank, the precise rank of $G$ as well as the ranks and
  dimensions of the Sylow subgroups of $G$ can be recognized by a
  single sentence in the basic first-order language of groups.
\end{abstract}

\maketitle


\section{Introduction}
In~\cite{NiSeTe21}, Nies, Segal and Tent carried out an investigation
of the model-theoretic concept of finite axiomatizability in the
context of profinite groups.  For instance, a profinite group $G$ is
finitely axiomatizable within a class $\mathcal{C}$ of profinite
groups, with respect to the first-order language
$\mathcal{L}_\mathrm{gp}$ of groups, if there is a sentence
$\psi_{G,\mathcal{C}}$ in $\mathcal{L}_\mathrm{gp}$ such that the
following holds: a profinite group $H$ in $\mathcal{C}$ is isomorphic
to $G$ if and only if $\psi_{G,\mathcal{C}}$ holds true in $H$, in
symbols $H\models \psi_{G,\mathcal{C}}$.  More generally, one takes
interest in whether specific properties or invariants of profinite
groups, again within a given class $\mathcal{C}$, can be detected
uniformly by a single sentence in $\mathcal{L}_\mathrm{gp}$.

Our main interest is in finitely generated profinite groups.  Nikolov
and Segal~\cite{NiSe07} established that such groups are strongly
complete; loosely speaking, this means that the topology of a finitely
generated profinite group is already predetermined by the abstract
group structure.  Jarden and Lubotzky~\cite{JaLu08} used Nikolov and
Segal's finite width results for certain words to prove that every
finitely generated profinite group is `first-order rigid', i.e.,
determined up to isomorphism by its first-order theory, within the
class of profinite groups.  By restricting to \emph{finite}
axiomatizability, we probe for more delicate first-order properties
within suitable classes of finitely generated profinite groups.

In this paper we focus on the class of profinite groups of finite
Pr\"ufer rank, from now on `rank' for short.  This invariant is
connected to, but not to be confused with the minimal number of
generators: the rank of a profinite group $G$ is defined as
\[
  \rk(G) = \sup \{ \dd(H) \mid H \le_\mathrm{o} G \} = \sup \{ \dd(H)
  \mid H \le_\mathrm{c} G \} ,
\]
where $\dd(H)$ denotes the minimal number of generators of a
topological group $H$ and, as indicated, $H$ runs over all open or all
closed subgroups of~$G$.  It is not difficult to see that the rank of
$G$ is the supremum of the ranks of its finite continuous quotients,
i.e.,
$\rk(G) = \sup \{ \rk(G / N) \mid N \trianglelefteq_\mathrm{o} G \}$.
The rank plays a central role in the structure theory of $p$-adic Lie
groups.  It is known that, for any given prime $p$, a profinite group
$G$ admits the structure of a $p$-adic analytic group if and only if
$G$ is virtually a pro-$p$ group of finite rank.  The dimension
$\dim G$ of a $p$-adic analytic profinite group $G$ is the analytic
dimension of $G$ as a $p$-adic manifold; in fact, $\dim G \le \rk(G)$
and $\dim G$ coincides with the rank $\rk(U)$ of any uniformly
powerful open pro-$p$ subgroup $U$ of~$G$.  Further details and
related results about $p$-adic analytic pro-$p$ groups can be found
in~\cite{DiDSMaSe99}; the concise introduction~\cite{Kl11} summarizes
key aspects of the theory.

\smallskip

Loosely speaking, our aim is to show that, for every finite set of
primes~$\pi$, the precise rank $r$ as well as the ranks
$\mathbf{r} = (r_p)_{p \in \pi}$ and dimensions
$\mathbf{d} = (d_p)_{p \in \pi}$ of the Sylow pro-$p$ subgroups of any
pro-$\pi$ group $G$ of finite rank can be recognized by a single
sentence $\sigma_{\pi,r,\mathbf{r},\mathbf{d}}$ in the first-order
language of groups~$\mathcal{L}_\mathrm{gp}$.  The starting point for
our investigation is Proposition~5.1 in \cite{NiSeTe21} which states:
Given $r \in \N$, there is an $\mathcal{L}_\mathrm{gp}$-sentence
$\rho_{p,r}$ such that for every pro-$p$ group $G$, the following
implications hold
\[
  \rk(G) \leq r \quad \Rightarrow \quad G \models \rho_{p,r} \quad
  \Rightarrow \quad \rk(G) \le r(2+\log_{2}(r)).
\]

Our first theorem both strengthens and generalizes this result.  The
$p$-rank $\rk_p(G)$ of a profinite group $G$ is the common rank of all
Sylow pro-$p$ subgroups of~$G$.  A sentence $\phi$ in
$\mathcal{L}_\mathrm{gp}$ is called an
$\exists \forall \exists$-sentence if it results from a
quantifier-free formula $\phi_0$ by means of a sequence of
existential, universal and existential quantifications (in this
order), rendering the free variables of $\phi_0$ to be bound
in~$\phi$; compare with Example~\ref{exm:recognize-finite-quotient}.

\begin{teor} \label{thm:rank-intro} Let $\pi$ be a finite set of
  primes.  Let $r \in \N$ and let
  $\mathbf{r} = (r_p)_{p \in \pi}$ be a tuple in
  $\{0, 1, \ldots,r\}$.  Then there exists an
  $\exists \forall \exists$-sentence $\varrho_{\pi, r, \mathbf{r}}$ in
  $\mathcal{L}_{\mathrm{gp}}$ such that, for every pro-$\pi$ group
  $G$, the following are equivalent:
  \begin{enumerate}
  \item[(i)] $\rk(G) = r$, and $\rk_p(G) = r_p$ for every
    $p \in \pi$.
  \item[(ii)] $\varrho_{\pi,r,\mathbf{r}}$ holds in $G$, i.e.,
    $G \models \varrho_{\pi, r,\mathbf{r}}$.
    \end{enumerate}
\end{teor}

It is no coincidence that the sentences $\varrho_{\pi,r,\mathbf{r}}$
which we manufacture to prove the theorem depend on the given set of
primes~$\pi$.  A standard ultraproduct construction reveals that, for
every \emph{infinite} set of primes $\tilde \pi$ and $r \in \N$, there
is no $\mathcal{L}_\mathrm{gp}$-sentence $\vartheta_{\tilde\pi,r}$
which could identify, uniformly across $p \in \tilde\pi$, among
pro-$p$ groups $G$ those with rank $\rk(G) = r$; see
Proposition~\ref{prop:ultra-prod}.

In addition to Theorem~\ref{thm:rank-intro} we establish a
corresponding theorem which concerns the dimensions of the Sylow
subgroups of a profinite group of finite rank.

\begin{teor} \label{thm:dim-intro} Let $\pi$ be a finite set of
  primes.  Let $r \in \N$ and let $\mathbf{d} = (d_p)_{p \in \pi}$ be
  a tuple in $\{0, 1, \ldots,r\}$.  Then there exists an
  $\exists \forall \exists$-sentence $\tau_{\pi,r,\mathbf{d}}$ in
  $\mathcal{L}_\mathrm{gp}$ such that, for every pro-$\pi$ group $G$
  with $\rk(G)=r$, the following are equivalent:
  \begin{enumerate}
  \item[(i)] For every $p \in \pi$, the Sylow pro-$p$
    subgroups of $G$ have dimension $d_p$.
  \item[(ii)] $\tau_{\pi,r,\mathbf{d}}$ holds in $G$, i.e.,
    $G \models \tau_{\pi, r,\mathbf{d}}$.
  \end{enumerate}
\end{teor}

In combination, the two theorems provide the first-order sentences
$\sigma_{\pi,r,\mathbf{r},\mathbf{d}}$ with the properties promised
above.  It is remarkable that such sentences exist in the basic
language $\mathcal{L}_\mathrm{gp}$ of groups.  In connection with
$p$-adic analytic profinite groups, it is often necessary to employ
suitably expanded languages in order to capture part of the
topological or analytic structure; compare with~\cite{MaTe16}.  We do
not need to enlarge the language at all.  Moreover, the complexity of
$\sigma_{\pi,r,\mathbf{r},\mathbf{d}}$ remains within three
alterations of $\exists$- and $\forall$-quantifiers, even though the
sentences that we manufacture depend strongly on the given set of
primes~$\pi$.

As we will show, the proofs of Theorems~\ref{thm:rank-intro} and
\ref{thm:dim-intro} reduce, in a certain sense, to the simpler setting
of pronilpotent pro-$\pi$ groups, termed $\mathcal{C}_\pi$-groups by
Nies, Segal and Tent~\cite[Section~5]{NiSeTe21}.  We recall that, even
in the pronilpotent case, Sylow subgroups are not in general definable
and there is no standard reduction to pro-$p$ groups; this can be seen
from relative quantifier elimination results (down to positive
primitive formulas) for modules over rings.  Part of our task is to
develop appropriate tools to by-pass this obstacle.

\smallskip

Key to our approach for proving Theorems~\ref{thm:rank-intro} and
\ref{thm:dim-intro} are purely group-theoretic considerations leading
to Theorem~\ref{thm:rank-in-finite-quotient} and its corollary, about
profinite groups which are virtually pronilpotent and of finite rank.
Specialising to the setting of finite nilpotent groups, we can
formulate the central insight as follows.

\begin{teor} \label{thm:nilp-gp-intro}
  Let $G$ be a finite nilpotent group of rank $r = \rk(G)$.  Then
  \[
    \rk(G) = \rk \!\big( G/\Phi^{j(r)}(G) \big) \qquad \text{for}
    \quad j(r) = 2r+\lceil \log_2(r) \rceil +2,
  \]
  where $\Phi^{j(r)}(G)$ denotes the $j(r)$th iterated Frattini
  subgroup of~$G$.
\end{teor}

It is an open problem to identify, if at all possible, even smaller
canonical quotients which witness the full rank of a finite nilpotent
group.

Following a suggestion of Gonz\'{a}lez-S\'{a}nchez, we derive from a
result of H{\'e}thelyi and L{\'e}vai~\cite{HeLe03} a new description
of the dimension of a finitely generated powerful pro\nobreakdash-$p$
group; this is useful for establishing Theorem~\ref{thm:dim-intro},
but also of independent interest.

\begin{teor} \label{thm:powerful-intro} Let $G$ be a finitely
  generated powerful pro-$p$ group with torsion subgroup $T$, and let
  $\Omega_{\{1\}}(G) = \{ g \in G \mid g^p = 1 \}$ denote the set of
  all elements of order $1$ or $p$ in~$G$.  Then
  \[
    \dim(G) = \dd(G) - \log_p \lvert \Omega_{\{1\}}(G)
      \rvert = \dd(G) - \dd(T).
  \]
\end{teor}

\smallskip

With a view toward possible future investigations, we add a final
comment and a question.  Naturally one wonders whether `being of
finite rank' per se can be captured by a suitable first-order
sentence.  Results of Feferman and Vaught~\cite{FeVa59} imply that,
even for a fixed prime~$p$, there is no set $\Sigma_p$ of
$\mathcal{L}_\mathrm{gp}$-sentences (and in particular no single
sentence) which identifies among the collection of all pro-$p$ groups
those that possess finite rank.  Indeed, the class of pro-$p$ groups
of finite rank is closed under taking finite cartesian products, but
an infinite cartesian product of non-trivial pro-$p$ groups of finite
rank is not even finitely generated.  Therefore
\cite[Corollary~6.7]{FeVa59} shows that no $\Sigma_p$ with the desired
property exists.  However, a modified question suggests itself.  Given
$d \ge 2$, is there a set $\Sigma_{p,d}$ of
$\mathcal{L}_\mathrm{gp}$-sentences (possibly a single sentence) such
that the following holds for pro-$p$ groups $G$ with $\dd(G) \le d$:
the group $G$ has finite rank if and only if $G$ satisfies
$\Sigma_{p,d}$?

\begin{oss-tech}
  Our proofs for Theorems~\ref{thm:rank-intro} and \ref{thm:dim-intro}
  involve results of Lucchini~\cite{Lu97} and an observation
  of~Mazurov~\cite{Ma94} which currently rely on the classification of
  finite simple groups.  However, in suitable circumstances, e.g., if
  we restrict attention to prosoluble groups, the required ingredients
  are know to hold without use of the classification; compare with
  \cite[Section~5]{Lu89}.  If $2 \not\in \pi$, the Odd Order Theorem
  guarantees that all pro-$\pi$ groups are prosoluble.
\end{oss-tech}

\begin{org-not}
  In Section~\ref{sec:detect-rk} we prove
  Theorem~\ref{thm:rank-in-finite-quotient} and its corollary, which
  specialize to Theorem~\ref{thm:nilp-gp-intro}. In
  Example~\ref{exa:what-goes-wrong} we discuss limitations of our
  strategy; Proposition~\ref{prop:ultra-prod} shows that
  Theorem~\ref{thm:rank-in-finite-quotient} does not generalize to
  groups involving infinitely many primes.  In Section~\ref{sec:rank}
  we establish Theorem~\ref{thm:rank-intro}, in
  Section~\ref{sec:dimension} we prove
  Theorem~\ref{thm:powerful-intro} and deduce
  Theorem~\ref{thm:dim-intro}.

  Our notation is mostly standard and in line with current practice.
  For instance, $\Zent(G)$ denotes the centre of a group~$G$, and
  $C_n$ denotes a cyclic group of order~$n$.  The meaning of possibly
  less familiar terms, such as $\Phi(G)$ for the Frattini subgroup and
  $\Phi_p(G)$ for the $p$-Frattini subgroup of a group~$G$, are
  explained at their first occurence.  We deal exclusively with
  profinite groups.  Accordingly, notions such as the Frattini
  subgroup, the commutator subgroup or the subgroup generated by a
  given set are tacitly understood in the topological sense: in each
  case we mean the topological closure of the corresponding abstract
  subgroup.  Basic model-theoretic concepts which are employed without
  further reference are covered by standard texts such as~\cite{Ho93}.
\end{org-not}

\noindent \textit{Acknowledgements.}  The results form part of the
first author's PhD research project.  We thank Jon
Gonz\'{a}lez-S\'{a}nchez for drawing our attention to structural
results about finite powerful $p$-groups which led to
Theorem~\ref{thm:powerful-intro} and thus helped us to streamline our
proof of Theorem~\ref{thm:dim-intro}.  We thank Immanuel Halupczok for
valuable conversations on model-theoretic background; in particular,
these led us to formulate Proposition~\ref{prop:ultra-prod}.

\section{Detecting the rank in bounded quotients} \label{sec:detect-rk}

Every compact $p$-adic analytic group $G$ has finite rank and contains
an open normal powerful pro-$p$ subgroup~$F$.  Since $F$ is a pro-$p$
group, its Frattini subgroup $\Phi(F)$ coincides with $[F,F]F^p$ and
$F/\Phi(F)$ is elementary abelian.  Since $F$ is powerful, we know
that $\rk(F) = \dd(F) = \rk(F/\Phi(F))$;
see~\cite[Theorem~3.8]{DiDSMaSe99}.  Furthermore, the iterated
Frattini series $\Phi^j(F)$, $j \in \N$, of $F$ coincides with both
the lower $p$-series and the iterated $p$-power series of~$F$.  It
provides a base of neighbourhoods for $1$ in~$G$ consisting of open
normal subgroups.  Consequently, the rank of $G$ is given by
\[
  \rk(G) = \sup \{ \rk(G / \Phi^j(F)) \mid j \in \N \} = \max
  \{ \rk(G / \Phi^j(F)) \mid j \in \N \};
\]
in other words, $\rk(G)$ is the terminal value of the non-decreasing,
eventually constant sequence $\rk(G / \Phi^j(F))$, $j \in \N$.

It is natural to look for an upper bound for the smallest $j \in \N$
such that $\rk(G) = \rk(G/ \Phi^j(F))$, a bound that is, as far as
possible, independent of $p$ and any special features of the pair
$F \le G$.  Based on our current knowledge, the strongest possible
outcome could be that $\rk(G) = \rk(G/\Phi(F))$ holds without any
exceptions.  More modestly, one can ask for weaker bounds, possibly
contingent on additional information regarding $\rk(G)$.

We establish a result of the latter kind, which applies more generally
to profinite groups $G$ of finite rank that admit a pronilpotent open
normal subgroup~$F$.  We recall that the $p$-rank $\rk_p(G)$ of a
profinite group $G$ is simply the rank $\rk(P)$ of a Sylow pro-$p$
subgroup $P$ of~$G$.  Furthermore, we write $\Phi_p(G) = [G,G] G^p$
for the $p$-Frattini subgroup of~$G$; the $p$-Frattini quotient
$G/\Phi_p(G)$ is the largest elementary abelian pro-$p$ quotient of
the profinite group~$G$.

\begin{teor} \label{thm:rank-in-finite-quotient} Let $R \in\N$.
  Suppose that the profinite group $G$ has an open normal subgroup
  $F \trianglelefteq_\mathrm{o} G$ which is pronilpotent and such that
  each Sylow subgroup of $F$ is powerful.

  \smallskip
  
  \textup{(i)} If $\rk_p(G) \leq R$ for some prime $p$, then
  \[
    \rk_p(G) = \rk_p \!\big( G/\Phi^{2R+1}(F) \big).
  \]

  \textup{(ii)} If $\rk (G)\leq R$, then
  \[
    \rk(G) = \rk \!\big( G/\Phi^{2R+1}(F) \big).
  \]
\end{teor}

\begin{proof}
  It is convenient to write $F_i = \Phi^i(F)$ for $i \in \N$.

  \smallskip
  
  \noindent (i) Let $p$ be a prime such that $r_p = \rk_p(G) \leq R$.
  We show that $r_p = \rk_p(G/F_{2R+1})$.  Since $F$ is pronilpotent,
  its Hall pro-$p'$ subgroup $P'$ is normal in~$G$.  Working
  modulo~$P'$, we may assume without loss of generality that $F$ is a
  powerful pro-$p$ group.  In this situation $G$ is virtually a
  pro-$p$ group.  Clearly, $r_p \ge \rk_p(G/F_{2R+1})$ and, for a
  contradiction, we assume that $r_p > \rk_p(G/F_{2R+1})$.  Choose a
  pro-$p$ subgroup $H \le_\mathrm{o} G$ of minimal index among the
  open pro-$p$ subgroups of $G$ with $\dd(H) = r_p$.

  The sequence $\dd(HF_j/F_j)$, $j \in \N$, is non-decreasing
  and eventually constant, with final constant value~$\dd(H)$.  Since
  $\dd(H) = r_p < 2R+1$, we conclude that $\dd(HF_j/F_j)$,
  $j \in \N$, cannot be strictly increasing until it becomes
  constant. Hence there exists $j = j(H) \in \N$ such that
  \begin{multline} \label{equ:jk<} \dd(HF_j/F_j) =
    \dd(HF_{j+1}/F_{j+1}) < \dd(HF_{j+2}/F_{j+2}) \\ < \ldots <
    \dd(HF_{j+k+1}/F_{j+k+1}) = \dd(H)
  \end{multline}
  for suitable $k = k(H)$ with $1 \le k \le r_p \le R$.  In particular, this
  set-up implies that $j+k+1>2R+1$, hence $j > R$ and
  $2j \ge j+R+1 \ge j+k+1$. Consequently, we see that
  $[F_j,F_j] \subseteq F_{2j} \subseteq F_{j+k+1}$ and there is no
  harm in assuming that
  \begin{equation*}
    [F_j,F_j] = F_{2j} = 1.
  \end{equation*}
  This reduction renders $G$ finite, with abelian normal $p$-subgroups
  \[
    A = F_j \qquad \text{and} \qquad B = F_{j+1} = \Phi(F_j) = A^p.
  \]
  
  We set $l = \dd(H/(H \cap B)) = \dd(HB/B) < \dd(H) = r_p$ and choose
  generators $y_1, \ldots, y_l$ for $H$ modulo $H \cap B$ so that
  \[
    L = \langle y_1, \ldots, y_l \rangle \le
    H
  \]
  satisfies $L B = H B$.  Put $m = \dd(H) - l = r_p - l \ge 1$.  A
  collection of elements generates $H$ if and only if it generates the
  Frattini quotient $H/\Phi(H)$; the latter is elementary abelian,
  because $H$ is a $p$-group.  Thus the minimal generating set
  $y_1, \ldots, y_l$ modulo $H \cap B$ can be supplemented to a
  minimal generating set for $H$: there are $b_1, \ldots, b_m \in B$
  such that
  \[
    H = \langle y_1, \ldots, y_l, b_1, \ldots, b_m
    \rangle \qquad \text{with} \qquad \dd(H) = r_p = l + m.
  \]
  We put $M = \langle b_1,\ldots, b_m\rangle^{H} \trianglelefteq H$ so
  that $H = LM$.
  
  Choose $a_1, \ldots, a_m \in A$ with $b_i = a_i^{\, p}$ for
  $1 \le i \le m$ and set
  \[
    \widetilde{H} = \langle y_1, \ldots, y_l , a_1, \ldots, a_m
    \rangle \le G.
  \]
  We claim that $\widetilde{H}$ is a $p$-subgroup of $G$ such that
  \begin{equation}\label{equ:index-and-rank-H-tilde} 
    \lvert G : \widetilde{H} \rvert < \lvert G : H \rvert \qquad \text
    {and} \qquad \dd(\widetilde{H}) = r_p,
  \end{equation}
  which yields the required contradiction.

  \smallskip
 
  Clearly, $\widetilde{H} \le H A$ is a $p$-group and
  $H \subseteq \widetilde{H}$.  Moreover, we see that
  $HA = \widetilde{H} A = L A$.  We may assume without loss of
  generality that $G = L A$.  In this situation $G$ is a
  $p$\nobreakdash-group; furthermore, $L \cap A \trianglelefteq G$ is
  normal.  By construction, compare with~\eqref{equ:jk<}, we have
  $\dd(L/(L \cap A)) = \dd(HA/A) = \dd(HB/B) = l = \dd(L)$.  Thus
  $L \cap A \subseteq \Phi(L) \subseteq \Phi(H)$ and there is no harm
  in assuming $L \cap A = 1$.  This gives
  \[
    G = L \ltimes A, \qquad H = L \ltimes M \qquad \text{and} \qquad
    \widetilde{H} = L \ltimes \widetilde{M} \quad \text{for
      $\widetilde{M} = \langle a_1, \ldots, a_m
      \rangle^{\widetilde{H}}$.}
  \]

  We supplement $y_1, \ldots, y_l$ to a minimal generating set
  $y_1, \ldots, y_l , \widetilde{a}_1, \ldots, \widetilde{a}_n$ for
  the $p$-group $\widetilde{H}$, for suitable $n \in \{0,1,\ldots,m\}$
  and $\widetilde{a}_1, \ldots, \widetilde{a}_n \in \widetilde{M}$.
  The $p$-power map $g \mapsto g^p$ induces a surjective $L$-invariant
  homomorphism $\alpha \colon \widetilde{M} \to M$ between finite
  abelian $p$-groups.  This implies
  $\lvert \widetilde{M} \rvert > \lvert M \rvert$ and thus
  $\lvert G : \widetilde{H} \rvert < \lvert G : H \rvert$.
  Furthermore, using the identity map on $L$ in combination with
  $\alpha$, we obtain a surjective homomorphism from
  $\widetilde{H} = L \ltimes \widetilde{M}$ onto $L \ltimes M = H$.
  This shows that $r_p = \dd(H) \le \dd(\widetilde{H}) \le r_p$ and
  hence $\dd(\widetilde{H}) = r_p$, which
  completes the proof of~\eqref{equ:index-and-rank-H-tilde}.

  \medskip

  \noindent (ii) Now suppose that $\rk(G) \le R$.  Clearly, the
  maximal local rank
  \[
    \operatorname{mlr}(G) = \max \big( \{ \rk_p(G) \mid p\hspace{0.15cm} \text{prime} \}
    \big)
  \]
  is at most $\rk(G)$.  Conversely, Lucchini established
  in~\cite[Theorem~3 and Corollary~4]{Lu97} that
  \[
    \rk(G) \le \operatorname{mlr}(G)+1,
  \]
  with equality if and only if there are
  \begin{itemize}
  \item[$\circ$] an odd prime $p$ such
    that $r_p = \rk_p(G) = \operatorname{mlr}(G)$ and
  \item[$\circ$] an open subgroup $H \le_\mathrm{o} G$ and
    $N \trianglelefteq_\mathrm{o} H$ such that
    \[
      H/\Phi_p(N) \cong H/N \ltimes N/\Phi_p(N) \cong C_q \ltimes
      C_p^{\, \operatorname{mlr}(G)},
    \]
    where $H/N \cong C_q$ is cyclic of prime order $q \mid (p-1)$, the
    $p$-Frattini quotient
    $N/\Phi_p(N) \cong C_p^{\, \operatorname{mlr}(G)}$ is elementary
    abelian of rank $\operatorname{mlr}(G)$, and $H/N$ acts via
    conjugation faithfully on $N/\Phi_p(N)$ by power automorphisms
    (i.e., by non-zero homotheties if we regard $N/\Phi_p(N)$ as an
    $\F_p$-vector space).
  \end{itemize}
  For short let us refer within this proof to such a pair $(H,N)$ as a
  `runaway couple' for $G$ with respect to~$p$.

  By (i), we have
  $\operatorname{mlr}(G) = \operatorname{mlr}(G/F_{2R+1})$, and hence
  it suffices to show: if $G$ admits a runaway couple, then so does
  $G/F_{2R+1}$, in fact, with respect to the same prime.  Suppose that
  $(H,N)$ is a runaway couple for $G$ with respect to an odd prime~$p$
  so that $H/\Phi_p(N) \cong C_q \ltimes C_p^{\, r_p}$ as detailed
  above, with the additional property that $\lvert G : H \rvert$ is as
  small as possible.  Assume for a contradiction that $G/F_{2R+1}$
  does not admit a runaway couple.

  As in the proof of (i) there is no harm in factoring out the Hall
  pro-$p'$ subgroup $P'$ of~$F$, because $H \cap F \subseteq N$ and
  $H \cap P' \subseteq \Phi_p(N)$.  Consequently we may as well assume
  that $F \trianglelefteq_\mathrm{o} G$ is a powerful pro-$p$ group,
  which makes $G$ virtually a pro-$p$ group.

  As in the proof of (i), the sequence
  \[
    \dd \!\big( H/ \big( (H \cap F_j) \Phi_p(N) \big) \big) = \dd \!
    \big( HF_j/ \Phi_p(N)  F_j \big), \quad j \in \N,
  \]
  is non-decreasing and eventually constant, with final constant
  value
  \[
     \dd(H/\Phi_p(N)) = \dd(H) = r_p + 1< 2R+1.
  \]
  We use the same arguments as before to conclude that there exists
  $j = j(H)$ such that the analogue of~\eqref{equ:jk<} for
  $H/\Phi_p(N)$ holds and we reduce to the situation where
  $[F_j,F_j] = F_{2j} = 1$.  This reduction renders $G$ finite, with
  abelian normal $p$-subgroups
  \[
    A = F_j \qquad \text{and} \qquad B = F_{j+1} = \Phi(F_j) = A^p;
  \]
  furthermore, we have
  \begin{equation}\label{equ:l-equal-AB}
    l = \dd \!\big( N/ \big( (H \cap A) \Phi_p(N) \big) \big) =
    \dd\!\big( N/ \big( (H \cap B) \Phi_p(N) \big) \big) < \dd(N /
    \Phi_p(N) ) = r_p.
  \end{equation}
  It suffices to produce a runaway couple
  $(\widetilde{H},\widetilde{N})$ for the group $HA$ with respect to
  $p$ such that
  $\lvert HA : \widetilde{H} \rvert < \lvert HA : H \rvert$; thus we
  may assume that
  \[
    G = HA.
  \]
  
  This reduction allows us to conclude that
  $\Phi_p(N) \cap A \trianglelefteq G$ and there is no harm in
  assuming $\Phi_p(N) \cap A = 1$.  Likewise
  $M = H \cap A \trianglelefteq G$, and reduction modulo $\Phi_p(N)$
  induces an embedding of $M \le N$ into the elementary abelian group
  $N/\Phi_p(N) \cong C_p^{\, r_p}$.  Using~\eqref{equ:l-equal-AB}, we
  conclude that
  \[
    M = H \cap A = H \cap B = \langle b_1, \ldots, b_m \rangle \cong
    C_p^{\, m} \qquad \text{for $m = r_p - l \ge 1$.}
  \]
  The normal subgroup $M \Phi_p(N) \trianglelefteq H$ decomposes as a
  direct product $M \times \Phi_p(N)$.  Recall that
  $H/\Phi_p(N) \cong C_q \ltimes C_p^{\, r_p}$, with the action given
  by power automorphisms.  We build a minimal generating set
  $x,y_1, \ldots, y_l, b_1, \ldots, b_m$ for $H$ modulo $\Phi_p(N)$ by
  choosing
  \[
    x \in H \smallsetminus N \quad \text{and} \quad y_1, \ldots, y_l
    \in N
  \]
  which supplement $b_1, \ldots, b_m$ suitably.  We set
  \[
    L_1 = \langle x, y_1, \ldots, y_l \rangle \le H \qquad \text{and}
    \qquad L = L_1 \Phi_p(N) \le H.
  \]
  
  In this situation $H = LM$ and we claim that $L \cap M = 1$ so that
  \[
    H = L \ltimes M.
  \]
  Indeed, our construction yields that the intersection in
  $H/\Phi_p(N) \cong C_q \ltimes C_p^{\, l+m}$ of the subgroups
  \[
    L/\Phi_p(N) = \langle \overline{x} \rangle \ltimes \langle
    \overline{y_1}, \ldots, \overline{y_l} \rangle \cong C_q \ltimes
    C_p^{\, l} \quad \text{and} \quad M \Phi_p(N)/\Phi_p(N) \cong M
    \cong C_p^{\, m}
  \]
  is trivial.  This
  gives $L \cap M \subseteq \Phi_p(N)$ and consequently
  $L \cap M \subseteq \Phi_p(N) \cap M = 1$.
  
  Put
  $\widetilde{M} = \{ a \in A \mid a^p \in M \} \trianglelefteq G$.
  Recall that $M = H \cap B$ and $B = A^p$.  The $p$-power map
  constitutes a surjective $G$-equivariant homomorphism
  $\widetilde{M} \to M$ whose kernel $K \trianglelefteq G$, say,
  includes~$M$.  From $L \cap M = 1$ we conclude that
  $LK \cap \widetilde{M} = (L \cap \widetilde{M}) K \subseteq K$.
  Moreover, we have $L \cap K \subseteq H \cap A = M$ and thus
  $L \cap K \subseteq L \cap M = 1$.

  These considerations show that the group
  $\widetilde{H} = L \widetilde{M}$ maps onto
  \[
    \widetilde{H}/K \cong LK/K \ltimes \widetilde{M}/K \cong
    L \ltimes M = H,
  \]
  and hence onto $C_q \ltimes C_p^{\, r_p}$.  Thus $\widetilde{H}$
  gives rise to a runaway couple for $G$, with respect to the prime
  $p$, just as $H$ does.  To conclude the proof we observe that
  $\lvert K \rvert \ge \lvert M \rvert \ge p$ implies
  $\lvert \widetilde{H} \rvert > \lvert \widetilde{H} \rvert / \lvert
  K \rvert = \lvert H \rvert$ and hence
  $\lvert G : \widetilde{H} \rvert < \lvert G : H \rvert$.
\end{proof}

\begin{cor} \label{cor:rank-in-finite-quotient} Let $R \in\N$.
  Suppose that the profinite group $G$ has an open normal subgroup
  $F \trianglelefteq_\mathrm{o} G$ which is pronilpotent.

  \smallskip
  
  \textup{(i)} If $\rk_p(G) \leq R$ for some prime $p$, then
  \[
    \rk_p(G) = \rk_p \!\big( G/\Phi^{2R+\lceil \log_2(R) \rceil +2}(F)
    \big).
  \]

  \textup{(ii)} If $\rk (G)\leq R$, then
  \[
    \rk(G) = \rk \!\big( G/\Phi^{2R + \lceil \log_2(R) \rceil +2}(F) \big).
  \]
\end{cor}

\begin{proof}
  As in the proof of Theorem~\ref{thm:rank-in-finite-quotient}, one
  reduces to the case in which $F$ is a pro-$p$ group for a single
  prime~$p$.  From $\rk(F) \le R$ it follows that
  $\Phi^{\lceil \log_2(R)\rceil +1}(F) \trianglelefteq_\mathrm{o} G$
  is powerful; compare with~\cite[Chapter~2, Exercise~6]{DiDSMaSe99}.  Thus
  we can apply Theorem~\ref{thm:rank-in-finite-quotient} to
  $\Phi^{\lceil \log_2(R) \rceil +1}(F)$ in place of~$F$.
\end{proof}

The following example puts the basic idea behind the proof of
Theorem~\ref{thm:rank-in-finite-quotient} into perspective.  It
indicates that one would need to take a different approach or at least
make more careful choices in order to eliminate the dependency on the
parameter~$R$.

\begin{ex} \label{exa:what-goes-wrong} Let $n \in \N$ and
  consider the metabelian pro-$p$ group
  \[
    G = C \ltimes A, \qquad \text{where
      $C = \langle c \rangle \cong \Z_p$,} \quad
    \text{$A = \langle a_1, \ldots, a_n \rangle \cong \Z_p^{\,
        n}$}
  \]
  and the action of $C$ on $A$ is given by
  \[
    a_i^{\, c} = a_i a_{i+1} \quad \text{for $1 \le i < n$}, \qquad
    \text{and} \quad  a_n^{\, c} = a_n.
  \]
  Here $\Z_p$ denotes the additive group of the $p$-adic integers,
  viz.\ the infinite procyclic pro-$p$ group.  Then
  $G = \langle c, a_1 \rangle$ is $2$-generated, nilpotent of
  class~$n$ and has rank $\rk(G) = n+1$.  For instance,
  \[
    H = \langle c, a_1^{\, p^{n-1}}, a_2^{\, p^{n-2}}, \ldots,
    a_{n-1}^{\, p}, a_n \rangle \le_\mathrm{o} G
  \]
  requires $n+1$ generators.

  Suppose that $p > n \ge 2$.  Then
  $F = \langle c^p \rangle \ltimes A \trianglelefteq_\mathrm{o} G$ is
  powerful, and
  $\Phi^j(F) = \langle c^{p^j} \rangle \ltimes A^{p^{j-1}}$ for
  $j \in \N$.  Thus any subgroup
  $\widetilde{H} \le_\mathrm{o} G$ with
  $\widetilde{H} F = H F = \langle c \rangle F$ and
  $\dd(\widetilde{H}) = \dd(\widetilde{H} \Phi^n(F) / \Phi^n(F))$
  requires less than $\dd(H) = n+1$ generators, but nevertheless
  $\rk(G) = \rk(G/ \Phi(F))$.  The group
  \[
  K = \langle c^p, a_1 , \ldots, a_n \rangle,
  \]
  which is unrelated to~$H$, requires $n+1$ generators, even
  modulo~$\Phi(F)$.
\end{ex}

\section{Finite axiomatizability of the rank} \label{sec:rank}

In this section we establish Theorem~\ref{thm:rank-intro}.  We begin
with a basic example which illustrates the concept of an
$\exists \forall \exists$-sentence in $\mathcal{L}_\mathrm{gp}$ and
related constructions which we use frequently; compare
with~\cite[Sections~2 and~5]{NiSeTe21}.

\begin{ex} \label{exm:recognize-finite-quotient} Let $G$ be a
  profinite group and let $N \subseteq_\mathrm{c} G$.  Suppose that
  $N$ is definable in~$G$; this means that there is an
  $\mathcal{L}_\mathrm{gp}$-formula $\varphi(x)$, with a single free
  variable~$x$, such that $N = \{ g \in G \mid \varphi(g) \}$.

  Let $B = \{ b_1, \ldots, b_n \}$ be a finite group of order~$n$,
  with multiplication `table'
  \[
    b_i b_j = b_{m(i,j)}
  \]
  encoded by a suitable function
  $m \colon \{1,\ldots,n\} \times \{1,\ldots,n\} \to \{1,\ldots,n\}$.

  Then the $\exists \forall$-sentence
 \begin{multline*}
   \exists a_1, \ldots a_n \, \forall x,y, z : \; \varphi(1) \,\land\,
   \Big( \big( \varphi(x) \land \varphi(y) \big) \to \varphi( x^{-1} y
   ) \Big) \,\land\, \Big( \varphi(x) \to \varphi \big(y^{-1} x y \big)
   \Big) \\
   \land\, \Bigg( \bigwedge_{1 \le i < j \le n} \neg \varphi
   \big(a_i^{\, -1} a_j \big) \Bigg) \,\land\, \Bigg( \bigvee_{1 \le i \le
     n} \varphi \big(a_i^{\, -1} y \big) \Bigg) \,\land\, \Bigg(
   \bigwedge_{1 \le i,j \le n} \varphi \big( a_{m(i,j)}^{\, -1} a_i a_j
   \big) \Bigg)
  \end{multline*}
  can be used to certify that $N \trianglelefteq_\mathrm{c} G$ and
  $G/N \cong B$.  In particular, if $N \subseteq_\mathrm{c} G$ is
  $\exists$-definable, i.e., definable by means of an
  $\exists$-formula, we obtain an $\exists \forall \exists$-sentence
  to certify that $N \trianglelefteq_\mathrm{c} G$ and $G/N \cong B$.

  For instance, if we know or suspect that the commutator word has a
  certain finite width in $G$, we may consider the $\exists$-definable
  set
  \[
    N = \{ [x_1,y_1] \cdots [x_r,y_r] \mid x_1, y_1, \dots,
    x_r, y_r \in G \} \subseteq_\mathrm{c} G,
  \]
  for a given parameter $r \in \N$, and formulate an
  $\exists \forall \exists$-sentence in $\mathcal{L}_\mathrm{gp}$
  which certifies that, indeed, $N$ is equal to the entire commutator
  subgroup $[G,G]$ and that the abelianization $G/[G,G]$ is isomorphic
  to a given finite group.

  \medskip
  
  Sometimes we want to express, by means of an
  $\mathcal{L}_\mathrm{gp}$-sentence, extra features of a definable
  subgroup $H \le_\mathrm{c} G$.  This process typically involves
  quantification over elements of $H$ rather than~$G$ which, in
  general, may increase the quantifier complexity of the resulting
  sentences.  However, if $H = \{ g \in G \mid \varphi(g) \}$ is
  $\exists$-definable, where $\varphi(x)$ takes the form
  $\exists \underline{z} \colon\varphi_0(x,\underline{z})$ with
  $\varphi_0$ quantifier-free in free variables $x$ and
  $z_1,\ldots,z_m$, say, then $H$ is `quantifier-neutral' in the
  following sense.  First-order assertions about $H$ can be translated
  into assertions of the same quantifier complexity about~$G$, simply
  by expressing universal quantification over elements of $H$ as
  $\forall x, \underline{z}: ( \varphi_0(x,\underline{z}) \to \ldots)$
  and existential quantification over elements of $H$ as
  $\exists x, \underline{z}: ( \varphi_0(x,\underline{z}) \land
  \ldots)$.
\end{ex}

It is convenient to establish the assertions of
Theorem~\ref{thm:rank-intro} first for pronilpotent groups before
considering the general situation.

\begin{prop} \label{prop:sentence-rank} Let $\pi$ be a finite set of
  primes, let $r\in\N$ and let $\mathbf{r} = (r_p)_{p \in \pi}$ be a
  tuple in $\{0,1,\ldots,r\}$.  Then there exists an
  $\exists \forall \exists$-sentence $\omega_{\pi, r, \mathbf{r}}$ in
  $\mathcal{L}_{\mathrm{gp}}$ such that, for every pronilpotent
  pro-$\pi$ group $H$, the following are equivalent:
  \begin{enumerate}
  \item[(i)] $\rk(H) = r$, and $\rk_p(H) = r_p$ for every
    $p \in \pi$.
  \item[(ii)] $\omega_{\pi,r,\mathbf{r}}$ holds in $H$, i.e.,
    $H \models \omega_{\pi, r,\mathbf{r}}$.
  \end{enumerate}
\end{prop}

\begin{proof}
  We set $k = \lvert \pi \rvert$, write $\pi = \{ p_1, \ldots, p_k \}$
  and put $q = q(\pi) = p_1 \cdots p_k$.  As $H$ is pronilpotent, it
  is the direct product $H = \prod_{i=1}^k H_i$ of its Sylow pro-$p_i$
  subgroups~$H_i$.  We set $m = m(r) = \lceil \log_2(r) \rceil +1$.

  Similar to Example~\ref{exm:recognize-finite-quotient}, there is an
  $\exists \forall \exists$-sentence $\beta_1$ in
  $\mathcal{L}_\mathrm{gp}$ to certify that there are elements
  $a_1, \ldots, a_r$ in $H$ such that every element $h \in H$ can be
  written as $h = \prod_{j=1}^r a_j^{\, e_j} b$, for suitable choices
  for $e_j \in \{0,1, \ldots, q-1\}$ and
  \[
    b \in B(H) = \{ [x_1,y_1] \cdots [x_r,y_r] z^q \mid x_1, y_1,
    \dots, x_r, y_r, z \in H \} \subseteq_\mathrm{c} \Phi(H).
  \]
  Using \cite[Lemma~1.23]{DiDSMaSe99} as in
  \cite[Section~5]{NiSeTe21}, we see that $\beta_1$ holds for $H$ if
  and only if $\dd(H) = \dd(H/\Phi(H)) \le r$; moreover, in this case
  $\Phi(H) = B(H)$.  Consequently, the subgroup $\Phi(H)$ is
  $\exists$-definable in~$H$ and hence quantifier-neutral in the sense
  of Example~\ref{exm:recognize-finite-quotient}.  By recursion, there
  is an $\exists \forall \exists$-sentence $\beta_{m+1}$ such that
  $\beta_{m+1}$ holds for $H$ if and only if
  \begin{equation} \label{equ:frattini-quotients-bd}
    \rk(\Phi^j(H)/\Phi^{j+1}(H)) \le r \qquad \text{for
      $0 \le j \le m$;}
  \end{equation}
  in this case the subgroup $F = \Phi^m(H)$ is $\exists$-definable
  in~$H$ and hence quantifier-neutral, moreover it satisfies
  $\dd(F) \le r$.  Furthermore, there is an $\forall \exists$-sentence
  $\gamma$ which certifies that $F$ is semi-powerful, in the
  terminology introduced in~\cite[Section~5]{NiSeTe21}: by
  \cite[Proposition~2.6]{DiDSMaSe99}, it suffices to express that
  every commutator $[x,y]$ of elements $x,y \in F$ is a $(2q)$th power
  $z^{2q}$ of a suitable $z \in F$.

  Once $F$ is $r$-generated and semi-powerful, we know that
  $\rk(F) \le r$.  If, in addition, the rank bounds specified
  in~\eqref{equ:frattini-quotients-bd} hold, we deduce that
  $\rk(H/F) \le m r$ and hence $\rk(H) \le R$ for $R = (m+1)r$.
  Furthermore, the group
  \[
    \Phi^{2R+1}(F) = \big\{ x^{q^{2R+1}} \mid x \in F \big\}
  \]
  is $\exists$-definable in~$H$ and hence quantifier-neutral; in
  particular, $H/\Phi^{2R+1}(F)$ is interpretable in~$H$.  Finally,
  $\lvert H/\Phi^{2R+1}(F) \rvert$ is bounded by $q^{(2R+m+1)r}$ and
  there is an $\exists \forall \exists$-sentence $\theta$ which
  certifies that $H/\Phi^{2R+1}(F)$ is one of the finitely many finite
  $\pi$-groups of suitable order which has rank~$r$ and whose
  $p$-ranks are in agreement with the prescribed~$\mathbf{r}$; compare
  with Example~\ref{exm:recognize-finite-quotient}.

  With the backing of Theorem~\ref{thm:rank-in-finite-quotient}, we
  form the conjunction of the sentences $\beta_{m+1}, \gamma, \theta$
  to arrive at an $\exists \forall \exists$-sentence
  $\omega_{\pi,r,\mathbf{r}}$ with the desired property.
\end{proof}

\begin{proof}[Proof of Theorem~\ref{thm:rank-intro}]
  We analyse the structure of a pro-$\pi$ group $G$ of rank
  $\rk(G) = r$ to build step-by-step a first-order sentence
  $\eta_{\pi,r}$ that is satisfied by any such group~$G$.  Following
  that we check that, conversely, every pro-$\pi$ group satisfying
  $\eta_{\pi,r}$ has rank at most~$2r$.  Applying
  Theorem~\ref{thm:rank-in-finite-quotient}, we extend $\eta_{\pi,r}$
  to a sentence $\varrho_{\pi,r,\mathbf{r}}$ which pins down precisely
  the rank as being $r$ and the ranks of the Sylow subgroups as being
  given by~$\mathbf{r}$.

  Our discussion involves upper bounds for certain integer parameters
  that depend on $\pi$ and $r$, but not on the specific group~$G$ used
  in our discussion; for short, we say that such parameters are
  $(\pi,r)$-bounded.

  \medskip
  
  \noindent \emph{Step 1.} The classification of finite simple groups
  implies that the set
  \begin{align*}
    \mathcal{S} = \mathcal{S}_{\pi,r} %
    & = \{ S\mid S \text{ a finite
      simple $\pi$-group such that $\rk(S) \le r$} \} \\
    & \subseteq \{
      S\mid S \text{ a finite simple $\pi$-group} \}
  \end{align*}
  is finite; see \cite[Remark following Lemma~2]{Ma94}.  Consequently,
  the cardinality of the set
  \[
    \Psi = \Psi_{G,\pi,r} =\{ \psi \mid
    \text{$\psi \colon G \to \Aut (S^l)$ a homomorphism for
      $S \in \mathcal{S}$ and $0 \le l \le r$} \}
  \]
  is $(\pi,r)$-bounded, because $G$ can be generated by at most $r$
  elements and any homomorphism between groups is determined by its
  effect on a chosen set of generators.  From this we observe that the
  index of
  \[
    K = K_{G,\pi,r} = \bigcap_{\psi \in \Psi} \Kern \psi
    \trianglelefteq_\mathrm{o} G
  \]
  in $G$ is $(\pi,r)$-bounded.  Thus there exists $f(\pi,r) \in \N$,
  depending on $\pi$ and $r$, but not on the specific group $G$, such
  that $\lvert G \colon K \rvert$ divides $f(\pi,r)$.

  We claim that $K$ is pronilpotent.  For this it suffices to show
  that $K/(K \cap L)$ is nilpotent for each
  $L \trianglelefteq_\mathrm{o} G$.  Let
  $L \trianglelefteq_\mathrm{o} G$.  By pulling back a chief series
  for the finite group $G/L$ to $G$, we obtain a normal series
  \[
    L = G_{n+1} \trianglelefteq G_n \trianglelefteq \ldots
    \trianglelefteq G_1 = G
  \]
  of finite length $n$ such that, for each $i \in \{1, \ldots, n \}$,
  the group $G_{i}/G_{i+1}$ is a minimal normal subgroup of
  $G/G_{i+1}$ and thus isomorphic to $S_i^{\, m(i)}$ for suitable
  choices of $S_i \in \mathcal{S}$ and $m(i) \in \N$.  Since each of
  the groups $S_i^{\, m(i)}$ contains an elementary abelian
  $p$-subgroup of rank~$m(i)$, for primes $p$ dividing
  $\lvert S_i \rvert$, we obtain
  $m(i) \le \rk \!\big( S_i^{\, m(i)} \big) \le \rk(G) = r$ for all
  $i \in \{1, \ldots, n\}$.  Intersecting with $K$, we obtain a series
  \begin{equation} \label{equ:series-for-K}
    K \cap L = K \cap G_{n+1} \trianglelefteq K \cap G_n
    \trianglelefteq \ldots \trianglelefteq K \cap G_1 = K
  \end{equation}
  consisting of $G$-invariant subgroups with factors
  $(K \cap G_i) / (K\cap G_{i+1}) \cong S_i^{\, l(i)}$ satisfying
  $0 \le l(i) \le m(i) \le r$, for $i \in \{1, \ldots, n\}$.  By
  construction, $K$ acts trivially on each of these factors so that
  $[K \cap G_i, K] \subseteq K \cap G_{i+1}$ for
  $i \in \{1, \ldots, n\}$.  Thus \eqref{equ:series-for-K} constitutes
  a central series for $K/(K \cap L)$, and $K/(K \cap L)$ is nilpotent
  (of class at most $n$).

  \medskip
  
  \noindent \emph{Step 2.} Next we consider the group
  \[
    H = G^{f(\pi,r)} = \langle g^{f(\pi,r)} \mid g \in G \rangle
    \trianglelefteq_\mathrm{o} G \qquad \text{with} \quad H \subseteq
    K;
  \]
  the index $\lvert G : H \rvert$ is $(\pi,r)$-bounded, by the
  positive solution to the Restricted Burnside Problem.  In fact, we
  do not require the general result, but a rather special case, which
  is easy to establish.  Indeed, assume for the moment that the
  pro-$\pi$ group $G$ of rank $r$ is finite of exponent $f(\pi,r)$.
  We need to show that $\lvert G \rvert$ is $(\pi,r)$-bounded.  In
  Step~1 we established that $G$ has a nilpotent normal subgroup $K$
  of $(\pi,r)$-bounded index.  Thus there is no harm in assuming that
  $G=K$.  Furthermore, $K$ is a direct product of its Sylow
  $p$-subgroups, where $p$ ranges over the finite set $\pi$.  Hence we
  may even assume that $G$ is a $p$-group of rank at most~$r$, for
  some $p \in \pi$, and that $f(\pi,r)$ is a $p$-power, $p^e$ say.  In
  this situation, $G$ contains a powerful normal subgroup of
  $(p,r)$-bounded index (see~\cite[Theorem~2.13]{DiDSMaSe99}), and we
  may assume that $G$ itself is powerful.  The $p$-power series of a
  powerful $p$-group coincides with its lower $p$-series, and we
  obtain the bound $\lvert G \rvert \le p^{re}$.

  Next we observe that the verbal subgroup $H$ is an
  $\exists$-definable subgroup of~$G$ and hence quantifier-neutral, in the
  sense discussed in Example~\ref{exm:recognize-finite-quotient}.
  Indeed, by \cite[Theorem~1]{NiSe11}, every element of $H$ can be
  written as a product of a $(\pi,r)$-bounded number of $f(\pi,r)$th
  powers.  But again we only require the bound in a rather special
  case which is much easier to handle.  Indeed, descending without
  loss of generality to a subgroup of $(\pi,r)$-bounded index, as
  above, it suffices to recall that in a powerful pro-$p$ group every
  product of $p^e$th powers is itself a $p^e$th power; see
  \cite[Corollary~3.5]{DiDSMaSe99}.

  \medskip
  
  \noindent \emph{Step 3.}  Since $K$ is pronilpotent, so is $H$.  In
  the situation at hand, this fact can be expressed by an
  $\exists \forall \exists$-sentence.  Indeed, $H$ is pronilpotent if
  and only if $H / \Zent(H)$ is pronilpotent.  Hence it suffices to
  express the assertion that $H / \Zent(H)$ is pronilpotent.  Clearly,
  $\Zent(H)$ is $\forall$-definable in $H$ and hence in~$G$.  We set
  $k = \lvert \pi \rvert$ and write $\pi = \{ p_1, \ldots, p_k \}$.
  As $H$ is pronilpotent, $H = \prod_{i=1}^k H_i$ is the direct
  product of its Sylow pro-$p_i$ subgroups~$H_i$ and
  $\Zent(H) = \prod_{i=1}^k \Zent(H_i)$ so that
  $H / \Zent(H) \cong \prod_{i=1}^k H_i / \Zent(H_i)$.  From
  \[
    C_i = \Cen_H(H_i) = \prod\nolimits_{j=1}^{i-1} H_j \times \Zent(H_i)
    \times \prod\nolimits_{j=i+1}^k H_j, \qquad \text{for
      $i \in \{ 1, \ldots, k\}$,}
  \]
  we deduce that
  \[
    D_i = \bigcap \{ C_j \mid 1 \le j \le k \text{ and } j \ne i \} =
    \prod\nolimits_{j=1}^{i-1} \Zent(H_j) \times H_i \times
    \prod\nolimits_{j=i+1}^k \Zent(H_j)
  \]
  and thus
  \[
    D_i / \Zent(H) \cong H_i / \Zent(H_i), \qquad \text{for
      $i \in \{ 1, \ldots, k\}$.}
  \]
  As $\rk(G) \le r$, there exist, for each $i \in \{ 1, \ldots, k \}$,
  elements $x_{i,1}, \ldots, x_{i,r} \in H_i$ such that
  $H_i = \langle x_{i,1}, \ldots, x_{i,r} \rangle$ and thus
  \[
    C_i = \Cen_H(\{ x_{i,1}, \ldots, x_{i,r} \}).
  \]
  Subject to the $kr$ parameters $x_{1,1}, \ldots, x_{k,r}$, this
  makes $\Zent(H) = \bigcap_{i=1}^k C_i$ and each of the groups $D_i$
  quantifier-free definable, by suitable centralizer conditions;
  moreover $Q_i = D_i/\Zent(H)$ becomes interpretable in~$H$, for
  $1 \le i \le k$.

  We conclude that it suffices to express in an
  $\forall \exists$-sentence, subject to the $(\pi,r)$-bounded
  number of parameters $x_{s,t}$, that
  \begin{itemize}
  \item[(a)] $\bigcap_{i=1}^k C_i = \Zent(H)$, hence
    $\Zent(H) \subseteq D_i$, for $i \in \{1, \ldots k\}$;
  \item[(b)] $D_i/\Zent(H)$ is a pro-$p_i$ group for
    $i \in \{1, \ldots k\}$;
  \item[(c)] $[D_i,D_j] \subseteq \Zent(H)$ for $i,j \in \{1,\ldots,k\}$
    with $i \ne j$;
  \item[(d)]
    $H = D_1 \boldsymbol{\cdot} D_2 \boldsymbol{\cdot} \ldots
    \boldsymbol{\cdot} D_k$, where the right-hand side denotes the set
    of all products $y_1 \cdots y_k$ with factors $y_i \in D_i$ for
    $i \in \{1,\ldots,k\}$;
  \end{itemize}
  for this implies that $H / \Zent(H) = \prod_{i=1}^k D_i/\Zent(H)$ is
  the direct product of its Sylow subgroups and thus pronilpotent.
  Turning the parameters $x_{s,t}$ into variables bound by an extra
  existential quantifier at the front, we arrive at an
  $\exists \forall \exists$-sentence without parameters which verifies
  that $H$ is pronilpotent.

  Subject to the parameters $x_{s,t}$, the assertions in (a), (c) can
  be expressed by an $\forall$-sentence, and (d) can be achieved by
  means of an $\forall \exists$-sentence.  The only tricky part occurs
  in (b) where we need to express that the group $Q_i = D_i/\Zent(H)$
  is a pro-$p_i$ group.  Since we know a priori that $Q_i$ is a
  pro-$\pi$ group, this is achieved by demanding that every element of
  $Q_i$ is a $q_i$th power, for
  $q_i = p_1 \cdots p_{i-1} p_{i+1} \cdots p_k$.  This can be
  expressed by an $\forall \exists$-sentence at the level of~$H$,
  because $\Zent(H) = \bigcap_{i=1}^k C_i$ is quantifier-free
  definable subject to the parameters~$x_{s,t}$.
  
  \medskip

  \noindent \emph{Step 4.} By Step~2, the group $G / H$ is
  interpretable in $G$ and finite of $(\pi,r)$-bounded order.  There
  is an $\exists \forall \exists$-sentence that certifies that the
  factor group $G/H$ is among the finitely many finite groups of rank
  at most $r$ and exponent dividing $f(\pi,r)$; compare with
  Example~\ref{exm:recognize-finite-quotient}.  Using our results from
  Step~2, Step~3 and Proposition~\ref{prop:sentence-rank}, we produce
  an $\exists \forall \exists$-sentence that certifies that the power
  word $x^{f(\pi,r)}$ has $(\pi,r)$-bounded width in~$G$ and that
  $H = G^{f(\pi,r)}$ is pronilpotent of rank at most~$r$.

  The conjunction of these two sentences yields an
  $\exists \forall \exists$-sentence $\eta_{\pi, r}$ such that
  \begin{enumerate}
  \item[$\circ$] every pro-$\pi$ group $G$ of rank $\rk(G) = r$
    satisfies $\eta_{\pi, r}$;
  \item[$\circ$] conversely, if a pro-$\pi$ group $\tilde G$ satisfies
    $\eta_{\pi, r}$, then
    $\tilde H = \tilde G^{f(\pi,r)} \trianglelefteq_\mathrm{o} \tilde
    G$ is pronilpotent and both $\tilde H$ and $\tilde G / \tilde H$
    have rank at most $r$; in particular, this ensures that
    $\rk( \tilde G) \le R$ for $R = 2r$.
  \end{enumerate}

  We put $m = m(R) = \lceil \log_2(R) \rceil +1$.  As in the proof of
  Proposition~\ref{prop:sentence-rank} we see that
  $F = \Phi^{m(R)}(H) \trianglelefteq_\mathrm{o} G$ is
  $\exists$-definable, hence quantifier-neutral, and semi-powerful.
  Furthermore, $\Phi^{2R+1} (F)$ is $\exists$-definable, hence
  quantifier-neutral, and, by
  Theorem~\ref{thm:rank-in-finite-quotient},
  \[
    \rk(G) = \rk\! \big( G / \Phi^{2R+1} (F) \big) \quad \text{and}
    \quad     \rk_p(G) = \rk_p\! \big( G / \Phi^{2R+1} (F) \big)
    \text{ for every $p \in \pi$.}
  \]
  Just as in the proof of Proposition~\ref{prop:sentence-rank} we find
  an $\exists \forall \exists$-sentence which in conjunction with
  $\eta_{\pi, r}$ produces an $\exists \forall \exists$-sentence
  $\varrho_{\pi,r,\mathbf{r}}$ with the desired property.
\end{proof}

The next result complements Theorem~\ref{thm:rank-intro}.  It
illustrates that the rank of a pro-$p$ group cannot be detected by a
first-order sentence uniformly across all primes~$p$, even if the
language $\mathcal{L}_\mathrm{gp}$ was to be enlarged by an extra
function to be interpreted as the $p$-power map $x \mapsto x^p$ in
pro-$p$ groups.  We sketch a proof for completeness; it relies on a
standard ultraproduct construction and a well-known quantifier
elimination result in model theory.

\begin{prop}\label{prop:ultra-prod}
  Let $\tilde \pi$ be an infinite set of primes and let~$r \in \N$.
  Then there is no $\mathcal{L}_\mathrm{gp}$-sentence
  $\vartheta_{\tilde\pi,r}$ such that, for every $p \in \tilde\pi$
  and every finite elementary abelian $p$-group~$G$, the following are
  equivalent:
  \begin{enumerate}
  \item[(i)] $\rk(G) = r$.
  \item[(ii)] $\vartheta_{\tilde\pi,r}$ holds in $G$, i.e.,
    $G \models \vartheta_{\tilde\pi, r}$.
  \end{enumerate}
\end{prop}

\begin{proof}
  For a contradiction, assume that the
  $\mathcal{L}_\mathrm{gp}$-sentence
  $\vartheta = \vartheta_{\tilde\pi,r}$ has the desired property.
  Then $C_p^{\, r} \models \vartheta$ and
  $C_p^{\, r+1} \models \neg\vartheta$ for all $p \in \tilde\pi$.  We
  regard $C_p^{\, r}$ and $C_p^{\, r+1}$ as the additive groups of the
  vector spaces $\F_p^{\, r}$ and $\F_p^{\, r+1}$ over the prime
  field~$\F_p$.

  Let $\mathfrak{U}$ be a non-principal ultrafilter on the infinite
  index set $\tilde\pi$.  By {\L}o{\'s}'s theorem,
  \[
    \mathcal{K} = \left( \prod\nolimits_{p \in \tilde\pi} \F_p \right)
    / \sim_\mathfrak{U}
  \]
  is a field of characteristic~$0$, and
  \[
    \mathcal{V} = \left( \prod\nolimits_{p \in \tilde\pi} \F_p^{\, r}
    \right) / \sim_\mathfrak{U} \qquad \text{and} \qquad \mathcal{W} =
    \left( \prod\nolimits_{p \in \tilde\pi} \F_p^{\, r+1} \right) /
    \sim_\mathfrak{U}
  \]
  are non-zero $\mathcal{K}$-vector spaces.  Let
  $\mathcal{L}_{\mathcal{K}\text{-}\mathrm{vs}}$ denote the language
  of $\mathcal{K}$-vector spaces, which comprises the language of
  groups (for the additive group of vectors) and, for each
  scalar~$c \in \mathcal{K}$, a $1$-ary operation $f_c$ (to denote
  scalar multiplication by~$c$).  Clearly, the
  $\mathcal{L}_\mathrm{gp}$-sentence $\vartheta$ gives rise to an
  $\mathcal{L}_{\mathcal{K}\text{-}\mathrm{vs}}$-sentence~$\theta$,
  not involving scalar multiplication at all, such that by
  {\L}o{\'s}'s theorem
  \[
    \mathcal{V} \models \theta \qquad \text{and} \qquad \mathcal{W}
    \models \neg \theta,
  \]
  in contradiction to the known fact that the infinite
  $\mathcal{K}$-vector spaces $\mathcal{V}$ and $\mathcal{W}$ have the
  same theory, due to quantifier elimination; see
  \cite[Section~8.4]{Ho93}.
\end{proof}
  
\section{Finite axiomatizability of the
  dimension} \label{sec:dimension}

In this section we establish Theorems~\ref{thm:powerful-intro}
and~\ref{thm:dim-intro}.  We derive the former from a result of
H{\'e}thelyi and L{\'e}vai~\cite{HeLe03} about finite powerful
$p$-groups; compare with~\cite{Wi02,Fe07}.  We recall from
\cite[Theorem~4.20]{DiDSMaSe99} that the elements of finite order in a
finitely generated powerful pro-$p$ group form a powerful finite
subgroup, its torsion subgroup.

\begin{proof}[Proof of Theorem~\ref{thm:powerful-intro}]
  The torsion subgroup $T$ is finite and characteristic in $G$ so that
  $\Cen_G(T) \trianglelefteq_\mathrm{o} G$.  We choose a uniformly
  powerful open normal subgroup $U \trianglelefteq_\mathrm{o} G$ such
  that $U \subseteq \Cen_G(T)$ and $U \subseteq \Phi(G)$.  Since $U$
  is torsion-free, this implies that
  \[
    N = U \times T \trianglelefteq_\mathrm{o} G \qquad \text{and}
    \qquad \dd(G) = \dd(G/U).
  \]

  We show below that there exists $k \in \N$ such that
  $U^{p^k} = \Phi^k(U) \trianglelefteq_\mathrm{o} G$ satisfies
  \begin{equation} \label{equ:show-below}
    \Omega_{\{1\}} (G/U^{p^k}) = \Omega_{\{1\}} (N/U^{p^k}).
  \end{equation}
  Since $N/U^{p^k} \cong U/U^{p^k} \times T$ and because $U$ is
  uniformly powerful, $\Omega_{\{1\}} (N/U^{p^k})$ is in bijection
  with the cartesian product of sets
  \[
    \Omega_{\{1\}}(U/U^{p^k}) \times \Omega_{\{1\}}(T) =
    U^{p^{k-1}}/U^{p^k} \times \Omega_{\{1\}}(G)
  \]
  and furthermore
  $\log_p \lvert U^{p^{k-1}} / U^{p^k} \rvert = \dd(U)$.
  Put $s(G) = \log_p \lvert \Omega_{\{1\}}(G)
    \rvert$. Stringing all pieces together, we see that the finite
  powerful $p$-group $P = G/U^{p^k}$ satisfies
  \[
    \log_p \lvert \Omega_{\{1\}}(P) \rvert = \dd(U) + s(G) = \dim(G) +
    s(G).
  \]
  The theorem of H{\'e}thelyi and L{\'e}vai~\cite{HeLe03}
  yields
  $\log_p \lvert \Omega_{\{1\}}(P) \rvert = \dd(P)$ and
    $s(G) = \log_p \lvert \Omega_{\{1\}}(T) \rvert = \dd(T)$ so that
  \[
    \dim(G) = \log_p \lvert \Omega_{\{1\}}(P) \rvert - s(G) = \dd(P) -
    s(G) = \dd(G) - s(G) = \dd(G) - \dd(T).
  \]

  \smallskip

  It remains to establish~\eqref{equ:show-below}.  Since $U^{p^k}$,
  $k \in \N$, is a base for the neighbourhoods of $1$ in~$G$, it
  suffices to show that there exists an open normal subgroup
  $W \trianglelefteq_\mathrm{o} G$ such that for every
  $x \in G \smallsetminus N \subseteq_\mathrm{c} G$ we have
  $x^p \not\in W$, or in other words $x^p \not \equiv_W 1$.  From
  $T \subseteq N$ we see that $G \smallsetminus N$ does not contain
  any elements of finite order.  Hence for every
  $x \in G \smallsetminus N$ there exists
  $W_x \trianglelefteq_\mathrm{o} G$ such that
  $x^p \not\equiv_{W_x} 1$, and consequently $y^p \not\equiv_{W_x} 1$
  for all $y \in xW_x \subseteq_\mathrm{o} G$.  Since
  $G \smallsetminus N$ is compact, it is covered by a finite union of
  such cosets $xW_x$, i.e.,
  $G \smallsetminus N \subseteq \bigcup_{x \in X} xW_x$ with
  $\lvert X \rvert < \infty$.  This implies that
  $W = \cap_{x \in X} W_x \trianglelefteq_\mathrm{o} G$ has the
  required property.
\end{proof}

\begin{proof}[Proof of Theorem~\ref{thm:dim-intro}]
  Let $p \in \pi$ and put $d = d_p$.  It suffices to explain how one
  can build an $\exists \forall \exists$-sentence $\tau_{\pi,r,p,d}$
  in $\mathcal{L}_\mathrm{gp}$ which certifies that a
  pro\nobreakdash-$\pi$ group $G$ of rank $\rk(G) = r$ has Sylow
  pro-$p$ subgroup dimension~$d$.  As in the proof of
  Theorem~\ref{thm:rank-intro} we work with a general
  pro-$\pi$ group $G$ with $\rk(G) = r$ to concoct $\tau_{\pi,r,p,d}$.
  
  Using the same approach as in the proof of
  Theorem~\ref{thm:rank-intro}, we find an $\exists$-definable and
  hence quantifier-neutral subgroup $H \trianglelefteq_\mathrm{o} G$
  that is pronilpotent and has $(\pi,r)$-bounded index in $G$;
  moreover the arrangement can be certified by means of a suitable
  $\exists \forall \exists$-sentence.  We put
  $m = m(r) = \lceil \log_2(r) \rceil +1$.  In the proof of
  Proposition~\ref{prop:sentence-rank} we saw that we can use an
  $\exists \forall \exists$-sentence to describe that $\Phi^m(H)$ is
  semi-powerful and of $(\pi,r)$-bounded index in~$H$; in parallel we
  can realize $\Phi^m(H)$ as an $\exists$-definable and hence
  quantifier-neutral subgroup.  The Sylow subgroup dimensions do not
  change if we pass from $G$ to an open subgroup.  Replacing $G$
  by~$\Phi^m(H)$, we may therefore assume without loss of generality
  that $G$ itself is pronilpotent and semi-powerful.

  As $G$ is pronilpotent, $G$ is the direct product of its powerful
  Sylow subgroups; let $G_p$ denote the Sylow pro-$p$ subgroup and
  $T_p$ its torsion subgroup.  By Theorem~\ref{thm:powerful-intro} it
  suffices to produce an $\exists \forall \exists$-sentence which pins
  down within the finite range $\{0,1,\ldots,r\}$ the invariants
  \[
    \dd(G_p) = \log_p \lvert G_p : \Phi(G_p) \rvert \qquad \text{and}
    \qquad \dd(T_p) = \log_p \lvert \Omega_{\{1\}}(G_p) \rvert,
  \]
  where $\Omega_{\{1\}}(G_p) = \{ g \in G_p \mid g^p = 1 \}$ is the
  set of all elements of order $1$ or~$p$.  We observe that
  $G_p / \Phi(G_p) \cong G / \Phi_p(G)$ is essentially the
  $p$-Frattini quotient of $G$ and that
  $\Omega_{\{1\}}(G_p) = \{ g \in G \mid g^p = 1 \}$.
  
  The Frattini quotient $G/\Phi(G)$ has $(\pi,r)$-bounded order and
  maps onto the $p$-Frattini quotient $G/\Phi_p(G)$.  As in the proof
  of Proposition~\ref{prop:sentence-rank}, the group $G / \Phi(G)$ is
  interpretable in~$G$.  There is an
  $\exists \forall \exists$-sentence which detects any prescribed
  isomorphism type of $G/\Phi(G)$ among a $(\pi,r)$-bounded number of
  possibilities; compare with
  Example~\ref{exm:recognize-finite-quotient}.  Forming a suitable
  disjunction, we can also detect the isomorphism type of the
  $p$-Frattini quotient $G/\Phi_p(G)$ and hence the minimal numbers of
  generators~$\dd(G_p)$.

  Clearly, the closed subset
  $\{ g \in G \mid g^p = 1 \} \subseteq_\mathrm{c} G$ is
  quantifier-free definable in~$G$.  Moreover, its size
    equals $p^{\dd(T_p)}$ and is thus at most $p^r$.  We can easily
  identify by means of an $\exists \forall$-sentence its precise size
  and hence the invariant $\dd(T_p)$.
\end{proof}


\end{document}